\documentclass{article}

\usepackage{makeidx}         
\usepackage{graphicx}        
\usepackage{multicol}        
\usepackage[bottom]{footmisc}

\usepackage{amsfonts}
\usepackage{amsmath}
\usepackage{xspace}
\newtheorem{algorithm}{Algorithm}
\newtheorem{proposition}{Proposition}

\newcommand{\propref}[1]{Proposition~\ref{#1}}

\newcommand{\secref}[1]{Section~\ref{#1}}

\newcommand{\al}{\alpha}

\newcommand{\la}{\lambda}

\newcommand{\eps}{\varepsilon}

\newcommand{\Aa}{{\cal A}}
\newcommand{\Pp}{{\cal P}}
\newcommand{\Kk}{{\cal K}}
\newcommand{\Ll}{{\cal L}}

\newcommand{\Ss}{{\cal S}}
\newcommand{\RR}{\mathbb{R}}

\newcommand{\ra}{\rightarrow}

\newcommand{\psdp}{\proj_{\sdp}}
\newcommand{\pK}{\proj_{\Kk}}
\newcommand{\pA}{\proj_{\Aa}}

\newcommand{\trans}[1]{{#1}^{\top}}

\DeclareMathOperator{\proj}{Proj}
\DeclareMathOperator{\prox}{Prox}
\DeclareMathOperator{\trace}{trace}

\DeclareMathOperator{\rank}{rank}
\DeclareMathOperator{\Int}{int}
\DeclareMathOperator{\Diag}{Diag}
\DeclareMathOperator{\diag}{diag}

\DeclareMathOperator{\textor}{or}

\newcommand{\sedumi}{SeDuMi\xspace}
\newcommand{\mprw}{MPRW\xspace}
\newcommand{\sdpnal}{SDPNAL\xspace}

\def\argmin{\mathop{\rm argmin}}

\newcommand{\val}{\ensuremath{\mathop{\rm val}}} 

\newcommand{\qqandqq}{\qquad\text{and}\qquad}

\newcommand{\Sy}{{\cal S}_n}
\newcommand{\sdp}{{\cal S}_n^{+}}

\renewcommand{\AE}{A_{\mathrm{E}}}
\newcommand{\AI}{A_{\mathrm{I}}}
\newcommand{\bE}{b_{\mathrm{E}}}
\newcommand{\bI}{b_{\mathrm{I}}}
\newcommand{\mE}{{m_{\mathrm{E}}}}
\newcommand{\mI}{{m_{\mathrm{I}}}}

\newcommand{\accol}[1]{\left\{\begin{array}{l} #1 \end{array}\right.}

\newcommand{\norm}[1]{\ensuremath{\Arrowvert #1 \Arrowvert}} 

\newcommand{\prods}[2]{{\langle #1,#2 \rangle}} 

\title{Projection methods in conic optimization}

\begin{document}

\author{Didier Henrion$^{1,2}$ \and
J\'er\^ome Malick$^3$}

\footnotetext[1]{CNRS; LAAS; 7 avenue du colonel Roche, F-31077 Toulouse, France;
Universit\'e de Toulouse; UPS, INSA, INP, ISAE; LAAS; F-31077 Toulouse, France.}

\footnotetext[2]{Faculty of Electrical Engineering, Czech Technical University in Prague,
Technick\'a 4, CZ-16626 Prague, Czech Republic.}

\footnotetext[3]{CNRS; LJK; Grenoble, France.}

\date{}

\maketitle

\begin{abstract}
There exist efficient algorithms to project a point onto the intersection of a convex cone and an affine subspace. Those conic projections are in turn the work-horse of a range of algorithms in conic optimization, having a variety of applications in science, finance and engineering. This chapter reviews some of these algorithms, emphasizing the so-called regularization algorithms for linear conic optimization, and applications in polynomial optimization. This is a presentation of the material of several recent research articles; we aim here at clarifying the ideas, presenting them in a general framework, and pointing out important techniques.
\end{abstract}

\section{Introduction, motivations, outline}

\subsubsection{Projection onto semidefinite positive matrices}
Consider the space $\Sy$ of symmetric $n$-by-$n$ matrices, equipped
with the norm associated to the usual inner product 
\[
\prods{X}{Y}=\sum^n_{i,j=1}X_{ij}Y_{ij}=\trace(\trans{X}Y).
\]
The subset $\sdp$ made of positive semidefinite matrices forms a
closed convex cone of $\Sy$. A general result for closed convex sets
yields that we can project onto~$\sdp$: given $C\in \Sy$, there exists
an unique element of $\sdp$ (called the projection of $C$ onto $\sdp$
and denoted by $\psdp(C)$) such that
\[
\norm{\psdp(C)-C} = \min_{X\in \sdp}\norm{X-C}.
\]
It turns out that we also have an explicit expression
of this projection, through the spectral decomposition of $C$.
Consider indeed the decomposition
\[
C = U \Diag(\la_1 ,\ldots, \la_n )\trans{U}
\] 
where  $\la_1 \geq \cdots \geq\la_n$ are the eigenvalues of $C$ 
and $U$ is a cor\-res\-pon\-ding orthonormal matrix of eigenvectors of $C$;
then the projection of $C$ onto $\sdp$ is 
\begin{equation}\label{projsdp}
\psdp(C) = U \Diag\big(\max(0, \la_1 ),\ldots,\max(0, \la_n )\,\big)\trans{U}.
\end{equation}
This result was noticed early by statisticians \cite{stat-1979} (see
also \cite{higham-1988}), and since then this projection has been
widely used. We notice that this result generalizes nicely for
``spectral sets''; see \cite{lewis-malick-2007}. Note also that the
numerical cost of computing this projection is essentially that of
computing the spectral decomposition of $C$, the matrix to project.

The developments of this chapter show that more sophisticated
projections onto subsets of $\sdp$ are also computable using
standard tools of numerical optimization. More specifically, the
subsets that we consider are intersections of the cone $\sdp$ with a
polyhedron (defined as affine equalities and inequalities). Though the
projection onto those intersections is not explicit anymore, we still
have efficient algorithms to compute them, even for large-scale
problems.

\subsubsection{Projection onto correlation matrices}
The most famous example of such projections is the projection onto the
set of correlation matrices (that are the real symmetric positive
semidefinite matrices with ones on the diagonal). It is common to be
faced with a matrix that is supposed to be a correlation matrix but
for a variety of reasons is not. For example, estimating correlation
matrices when data come from different time frequencies may lead to a
non-positive semidefinite matrix. Another example is stress-testing in
finance: a practitioner may wish to explore the effect on a portfolio
of assigning certain correlations differently from the historical
estimates, but this operation can destroy the semidefiniteness of the
matrix.

These important practical questions have led to much interest in the
problem of computing the nearest correlation matrix to a given a
matrix $C$ (see
e.g.\,\cite{higham-2002}, \cite{malick-2004}, \cite{qi-sun-2006} and
\cite{borsdorf-higham-2008}). This problem is simply formulated as the
projection of $C$ onto correlation matrices
\begin{equation}\label{eq:corr}
\accol{
    \min \quad  \frac{1}{2}\norm{X-C}^2\\
    \quad X_{ii} = 1, \quad i=1,\ldots,n\\
    \quad X\succeq 0.
}
\end{equation}
The methods reviewed in this chapter apply to solving this problem in
particular. The point is that this problem (and variants of it)
can now be solved efficiently (for sizes up to $n=5000$; the only
limitation on a standard computer is the memory constraints).

\subsubsection{Conic projection problem}
The general problem that we first focus on in this chapter is the following.
In the space $\RR^n$ equipped with the standard inner product, we want to
compute the projection of a point $c\in \RR^n$ onto the intersection
$\Kk\cap \Pp$ where
\begin{itemize}
\item $\Kk$ is a closed convex cone in $\RR^n$ (that we further assume
  to have full dimension in $\RR^n$; that is, $\Int\Kk \neq \emptyset$),
\item $\Pp$ is a convex polyhedron defined by affine (in)equalities
  \[
  \Pp:=\left\{x\in \RR^n\!:\,\trans{a_i\!}x= (\textor \leq)\;b_i,
  \quad i=1,\ldots,m\right\}.
  \]
\end{itemize}
We suppose moreover that the intersection $\Kk\cap\Pp$ is nonempty, so
that the projection onto the closed convex set $\Kk\cap\Pp$ exists and
is unique (see e.g.\;\cite{hull-1993}).

The fact that $\Pp$ is defined by both equalities and inequalities
does not really matter in our developments. To simplify presentation,
one may take only equalities, so that $\Pp$ is an affine subspace. We
prefer to keep the above loose notation with both equalities and
inequalities, because it is closer to projection problems arising in
practice, and because it does not impact the basics of projection
algorithms. Adding positive slack variables for the
inequalities allows us to reformulate the problem as a projection onto
the intersection of an affine space with a cone of the form $\Kk\times
(\RR_+)^\mI$.

We note that in general one can project onto a polyhedron $\Pp$. For
the case when there are only (independent) equalities in the
definition (i.e.\;if $\Pp=\Aa$ is an affine subspace of the equation $Ax=b$
with a full-rank matrix $A$), we have the explicit expression of the
projection of $x$
\begin{equation}\label{eq-projU}
  \proj_{\Aa}(x)=x-\trans{A}[A\trans{A}]^{-1}(A x-b).
\end{equation}
For a general polyhedron $\Pp$, we still can compute the projection
$\proj_{\Pp}(x)$ efficiently using quadratic programming solvers (see
e.g.\,\cite{nocedal-wright-1999} or \cite{bgls-2003}).  In this
chapter, we make the practical assumption that it is also easy to
project onto $\Kk$.  Recall from above that we have an easy-to-compute
expression of the projection for $\Kk =\sdp$; it turns out to be also
the case for the second-order cone (or Lorentz cone)
\[
\Ll_n:=\left\{x\in \RR^n:\ \norm{(x_1,\ldots,x_{n-1})}\leq
x_n\right\}.
\]
Though it is easy to project onto $\Pp$ and also onto $\Kk$ by
assumption, the projection onto the intersection $\Pp\cap\Kk$ can
still be challenging. The difficulty comes from the presence of both
(affine and conic) constraints at the same time. We will see in
\secref{sec:sdls} that many numerical methods to compute the
projection onto the intersection use combinations of projections onto
$\Pp$ and $\Kk$ separately.

The geometrical projection problem has an obvious analytical
formulation as a least-squares problem, namely minimizing the
(squared) norm subject to the conic constraints and the affine
constraints:
\begin{equation}\label{eq:pb}
\accol{\min \quad \frac{1}{2}\norm{x-c}^2\\
  \quad x\in \Pp\cap \Kk.}
\end{equation}
For example, an important subclass of such problems are
semidefinite least-squares problems (i.e.\;when $\Kk=\sdp$):
\begin{equation}\label{eq:sdls}
\accol{\min \quad \frac{1}{2}\norm{X-C}^2\\ 
  \quad \prods{A_i}{X}= (\textor \leq)\;b_i, \ \ i=1,\ldots,m\\ 
  \quad X\succeq 0,}
\end{equation} 
for $C,A_i\in \Sy$. For example, the nearest correlation matrix
problem \eqref{eq:corr} is an instance of this later class. Notice
finally that problems \eqref{eq:pb}~and~\eqref{eq:sdls}
coincide formally with $x\in \RR^{n^2}$ collecting the rows of $X\in
\Sy$.  This also explains the slight abuse of notation when writing
$\Kk=\sdp$. We use this ambiguity $x\leftrightarrow X$ in particular
in \secref{sec:poly} to ease presentation of relaxations of polynomial
optimization problems.

\subsubsection{Linear conic optimization problem}

The second development of this chapter is about a more
standard topic: solving linear conic optimization problems. With the
above notation, these problems can be expressed as
\begin{equation}\label{eq:lin}
\accol{\min \quad \trans{c}{x}\\
  \quad x\in \Pp\cap \Kk.}
\end{equation}
As presented in this handbook and as well as in the first handbook
\cite{wolkowicz-saigal-vandenberghe-2000}, linear conic programming
has been a very active field of research spurred by many applications
and by the development of efficient methods. 

In this chapter, we explain how conic projections can be used to
develop a new family of algorithms for solving linear conic problems.
In fact, the projection problem \eqref{eq:pb} can be written as a linear
problem of the form \eqref{eq:lin} and then can be solved using usual
conic programming solvers (we come back to this at the beginning of
\secref{sec:sdls} and we explain why it is not a good idea to do so).
However we will also show the other way around: \secref{sec:prox} explains
that one can also solve the linear conic problem \eqref{eq:lin} 
by solving projection problems~\eqref{eq:pb}, more precisely with a
succession of (truncated) projection-like problems. So-called
regularization methods are presented, discussed and illustrated on
solving semidefinite relaxations of combinatorial optimization and polynomial
optimization problems having many constraints.

\subsubsection{Polynomial optimization}

Over the last decade, semidefinite programming has been used in
polynomial optimization, namely for deciding whether a multivariate
real polynomial is nonnegative, or, more generally, to minimize a
polynomial on a semialgebraic set (a set described by polynomial
inequalities and equations). A hierarchy of embedded linear
semidefinite relaxations (of the form \eqref{eq:lin}) can
be constructed to generate a monotone sequence of bounds on the global
minimum of a polynomial optimization problem. Asymptotic convergence
of the sequence to the global minimum can be guaranteed under mild
assumptions, and numerical linear algebra can be used to detect global
optimality and extract global minimizers.  The theory is surveyed
in \cite{laurent-2009} and \cite{lasserre-2009}; the potential
applications are numerous (see e.g.\;in control theory
\cite{henrion-garulli-2005} or signal processing
\cite{dumitrescu-2007}).

Section\;\ref{sec:poly} reports numerical experiments showing that
regularization algorithms based on projections outperform classical
primal-dual interior-point algorithms for solving semidefinite
relaxations arising when deciding whether a polynomial is
nonnegative, and for globally minimizing a polynomial.

\subsubsection{Objectives and outline of this chapter}
This chapter focuses on projection problems that have a simple
geometric appeal as well as important applications in engineering. We
give references to some of these applications, and we give emphasis on
polynomial optimization.

The first goal of this chapter is to sketch the different approaches
to solve~\eqref{eq:pb}. \secref{sec:sdls} is devoted to this
review, with an emplasis on dual methods (in
\secref{sec:dual}). The bottomline is that as soon as we can easily
project onto $\Kk$ (we have in mind $\Ll_n$ and $\sdp$ as well as direct
products of these), we have efficient algorithms to project onto the
intersection $\Kk\cap\Pp$.

The second goal of this chapter is to explain how to use these conic
projections to build a family of ``regularization'' methods for linear
conic programming. The approach uses standard optimization techniques
(proximal algorithms and augmented Lagrangian methods) and has been
recently developed for the case $\Kk=\sdp$. \secref{sec:reg} presents
it in a general framework and underlines the role of conic
projections. The final section presents some numerical
experiments with regularization methods on polynomial optimization
problems, showing the interest of the approach in that context.

This chapter is meant to be an elementary presentation of parts of the
material of several papers; among those, our main references are
\cite{malick-2004}, \cite{qi-sun-2006},
\cite{malick-povh-rendl-wiegele-2009}, \cite{henrion-malick-2009},
\cite{zhao-sun-toh-2010} and \cite{nie-2009}.  We aim at clarifying
the ideas, presenting them in a general framework, unifying notation,
and most of all, pointing out what makes things work. To this purpose,
we have to elude some technical points; in particular, we discuss
algorithms, but we do not to give convergence results.
We try to give precise references throughout the text on these
lacking points.

\section{Conic projections: algorithms and applications}\label{sec:sdls}

This section reviews the methods for solving the conic projection
problem \eqref{eq:pb}, presenting them in chronological order. We
sketch the main ideas and give references; we do not get
into much details. Discussions about convergence issues
and numerical comparisons are beyond the scope of this section.

Beside interior-point methods, the basic idea of all the approaches is
to somehow separate the two constraint-sets $\Kk$ and $\Pp$ and to use
the projections onto them successively: this is obvious for
alternating projections and alternating directions methods; it is also
the case for dual methods (we focus on this latter method in
\secref{sec:dual}). The point is that we can solve the conic
projection problem \eqref{eq:pb} efficiently (by dual algorithms in
particular).

To simplify presentation, we stick here with the projection problem
\eqref{eq:pb}, but the approach generalizes in two directions.  First,
we could replace the cone $\Kk$ by any closed convex set: in this more
general case, the developments are similar, with slightly more
complicated expressions of dual objects (a related work is
\cite{micchellli-utretas-1988}). Second, we could consider problems with strongly
convex quadratic objective functions, such as
\begin{equation}\label{eq:gene}
  \accol{\min \quad \trans{(x-c)}Q(x-c) + \trans{d}{x}\\
    \quad x\in\Kk\cap \Pp}
  \end{equation}
with $Q$ positive definite. Such problems can be phrased as projection
problems with respect to $\norm{x}_Q = \sqrt{\trans{x}\!Q\,x}$ the
norm associated to~$Q$. The practical technical assumption is then
that one can project onto $\Kk$ with respect to $\norm{\cdot}_Q$
(which is not easy in general).

\subsection{Computing conic projections}\label{sec:algo}

\subsubsection{Using linear conic programming}
A tempting method to solve \eqref{eq:pb} is to cast this projection
problem as a usual linear conic programming problem, so that we can
use the powerful tools developed for this case. There are several
ways to do so; a simple one consists in pushing down the objective
function with an additional variable $t$: \eqref{eq:pb} is indeed
equivalent to linear conic program
\[
\accol{\min \quad t \\
  \quad x\in \Pp\\
  \quad x-c = z\\
  \quad (x,(z,t))\in \Kk\times \Ll_{n+1}}
\]
where the variable $z\in \RR^n$ is then introduced to express the
additional second-order cone constraint appearing in the
constraints. This problem can be readily given to usual conic solvers,
for example interior-points methods, like \sedumi
\cite{sturm-1999} or \texttt{SDPT3} \cite{tutuncu-toh-todd-2003} under
Matlab. Unfortunately, adding $(z,t)$ makes the computational cost and
memory space needed by a standard primal-dual interior-point method
increase, and numerical testing confirms that the method is not viable
in general (as mentioned e.g.\;in \cite{higham-2002},\cite{toh-2007}).

We note furthermore that the projection problem \eqref{eq:pb} is a
quadratic conic programming problem, hence a special case of nonlinear
conic optimization problems. We could solve \eqref{eq:pb} by
algorithms and software devoted to nonlinear conic optimization
problems such as the penalization method of \cite{kocvara-stingl-2003}.
However those methods would not use the special structure of
\eqref{eq:pb}, and as the above approach by linear conic programming,
they would be efficient only for small-size projection problems. The
projection problems are important enough to design algorithms
specifically to them, as presented in the sequel. Note that
we are not aware of a tailored penalization algorithm for~\eqref{eq:pb}.

\subsubsection{Alternating projections}
The alternating projection method is an intuitive algorithmic scheme
to find a point in the intersection of two sets: it consists in
projecting the initial point onto the first set, then projecting the
new point onto the second set, and then projecting again the new point
onto the first and keep on projecting alternatively. In other words,
it consists in repeating:
\begin{equation}\label{eq:alternating}
\accol{x_{k+1} = \proj_{\Kk}(y_k)\\
y_{k+1} = \proj_{\Pp}(x_{k+1})\\
}
\end{equation}
If the two sets have a ``regular'' intersection, this algorithm
converges linearly to a point in $\Pp\cap \Kk$ and we know the speed
of convergence (for two convex sets, see e.g.\;\cite{deutsch-2001};
for the general case, see the local result of
\cite{lewis-luke-malick-2008}).

We can modify this simple alternating projection scheme by adding a
correction step (called Dykstra's correction \cite{dykstra-1983}) at
each iteration \eqref{eq:alternating}
\begin{equation}\label{eq:dykstra}
\accol{x_{k+1} = \proj_{\Kk}(z_k)\\
y_{k+1} = \proj_{\Pp}(x_{k+1})\\
z_{k+1} = z_k-(x_{k+1}-y_{k+1}).
}
\end{equation}
This modification ensures the convergence of the sequence $(x_k)_k$ to
the projection $\proj_{\Kk\cap\Pp}(c)$ -- and not only to a point in
the intersection $\Kk\cap \Pp$. This approach was proposed by
\cite{higham-2002} for the nearest correlation matrix
problem~\eqref{eq:corr}. It generalizes to \eqref{eq:pb} since
it is easy to project onto $\Pp$ and we assume that it is the same for $\Kk$.
We will see that dual methods and alternating direction methods can be
interpreted as variants of this basic geometrical method.


\subsubsection{Dual methods}
The conic programming problem \eqref{eq:pb} looks more complicated
than a usual conic programming problem with linear function instead of
a norm as objective function. It turns out that the strong convexity
of the objective function provides nice properties to the dual problem
that can then be solved efficiently.
    
The dual approach was proposed for the conic
least-squares problem \eqref{eq:pb} in \cite{malick-2004}, later revisited by
\cite{boyd-xiao-2005} for the case of $\Kk=\sdp$, and then enhanced by
\cite{qi-sun-2006} and \cite{borsdorf-higham-2008} for the projection
onto correlation matrices. In the next section, we give more details
and more references about this approach.

\subsubsection{Interior points}
As a convex optimization problem, \eqref{eq:pb} can be attacked with
the interior-point machinery \cite{nemirovski-nesterov-1994}, assuming
that both the cone $\Kk$ and its polar cone 
\[
\Kk^o:=\left\{s\in \RR^n: \ \trans{s}x \leq 0 \ \text{for all }x\in \Kk\right\}
\]
are equipped with so-called
self-concordant barriers (as is the case for $\Ll_n, \sdp$).
The approach consists in solving perturbed optimality conditions of
\eqref{eq:pb}. 
As any projection problem, notice that the optimality condition is
\[
\bar x \in \Pp\cap \Kk, \quad 
\trans{(c- \bar x)}{(x- \bar x)}\leq 0, \quad \text{for all }x\in \Pp\cap \Kk.
\]
To write down the optimality conditions more concretely, let us
make explicit the affine constraints with the help of
$\AE\in\RR^{n\times \mE}$ and $\AI\in\RR^{n\times \mI}$ as
\begin{equation}\label{eq:pb2}
\accol{\min \quad \norm{x-c}^2\\
\quad \AE x=\bE, \ \AI x\leq \bI\\
\quad x\in \Kk.}
\end{equation}
Under a non-degeneracy assumption (e.g.\;Slater condition, see next
section), the optimality conditions of \eqref{eq:pb2} give
the complementarity system
\[
\accol{x-c+u+\trans{\AE}y+\trans{\AI}z = 0\\[0.5ex]
\AE x = \bE, \ y\in \RR^{\mE}\\[1ex]
\AI x \leq \bI, \ z\in \RR_+^\mI, \ \trans{z}(\AI x - \bI)=0\\[1ex]
x\in \Kk, \  u\in \Kk^o, \ \trans{u}{x}=0.
}
\]
Roughly speaking, an interior-point approach consists in
perturbing the complementary equations above and keeping other
equations satisfied. (We will see that the forthcoming dual approach
goes exactly the other way around.)  A first interior-point method is
proposed in \cite{takouda} for the nearest correlation matrix problem
\eqref{eq:corr}. Interior-point methods for general quadratic SDP are
introduced and tested on projection problems \eqref{eq:pb} in
\cite{tutuncu-toh-todd-2006} and \cite{toh-2007}.
    
\subsubsection{Alternating directions}

The alternating direction method is a standard method in variational
analysis (see e.g.\;\cite{gabay-1976}), going back to
\cite{douglas-rachford-1956}.  This method was proposed by
\cite{adm-2009} for solving the semidefinite projection problem
\eqref{eq:sdls} and by \cite{adm-2009-2} for more general
quadratically constrained quadratic SDP. The idea of the method is to
exploit the separable structure of the problem, as follows. Let us
duplicate the variables to write the equivalent problem
\begin{equation}\label{eq:pb3}
\accol{\min \quad \frac{1}{2}\norm{x-c}^2+ \frac{1}{2}\norm{y-c}^2\\
\quad x=y \\
\quad x\in \Kk,\quad  y\in \Pp.}
\end{equation}
The alternating direction method applied to \eqref{eq:pb3} gives the
following scheme: consider the augmented Lagrangian function
\[
L(x,y;z) = \frac{1}{2}\norm{x-c}^2+ \frac{1}{2}\norm{y-c}^2 
-\prods{z}{x-y} + \frac{\beta}{2}\norm{x-y}^2;
\]
the minimization of $L$ with respect to primal variables $(x,y)$ is
decomposed in two steps, so that an augmented Lagrangian iteration is
\[
\accol{x_{k+1} = \argmin_{x\in \Kk}L(x,y_k,z_k)\\
y_{k+1} = \argmin_{y\in \Pp}L(x_{k+1},y,z_k)\\
z_{k+1} = z_k-\beta(x_{k+1}-y_{k+1}).
}
\]
It is not difficult to prove that the two above minimizations boil down to
projections, more specifically
\[
x_{k+1}=\pK\!\Big(\frac{\beta y_k+z_k+c}{1+\beta}\Big),\quad
y_{k+1}=\pA\!\Big(\frac{\beta x_{k+1}-z_k+c}{1+\beta}\Big).
\]
Thus the approach alternates projections onto $\Pp$ and $\Kk$ to
compute the projection onto $\Kk\cap\Pp$; it can thus be seen
as a modification of the simple alternating projection scheme
\eqref{eq:alternating}, with the same flavour as Dykstra modification~\eqref{eq:dykstra}.

\subsection{More on dual approach}\label{sec:dual}

\subsubsection{Apply standard machinery}

Let us give more details about the dual approach for solving
\eqref{eq:pb}. Following \cite{malick-2004}, we apply the standard
mechanism of Lagrangian duality to this problem; we refer to
\cite[Ch.\,XII]{hull-1993} and \cite[Ch.\,5]{boyd-vandenberghe-2004}
for more on this mechanism in general.

Let us consider the more explicit form \eqref{eq:pb2}, and denote also
by $A:=[\AE;\AI]$ and $b:=[\bE;\bI]$ the concatenation of the affine
constraints.  We dualize affine constraints only: introduce the
Lagrangian, a function of primal variable $x\in \Kk$ and dual variable
$(y,z)\in \RR^\mE\!\times\RR_+^\mI$
\begin{equation}\label{eq-def-lagr}
  L(x;y,z) := \frac{1}{2}\norm{c - x}^2 - y^{\!\top}(\AE x-\bE) - z^{\!\top}(\AI x-\bI),
\end{equation}
and the corresponding concave dual function
\begin{equation}\label{eq-def-theta}
  \theta(y,z) :=  \min_{x \in \Kk} L(x;y,z),
\end{equation}
which is to be maximized. There is no more affine constraint in the
above minimum, and it is easy to prove (\cite[Th.3.1]{malick-2004})
that the problem corresponds to a projection onto $\Kk$: there
exists a unique point which reaches the above minimum, namely
\begin{equation}\label{eq-xy}
  x(y,z) :=\pK(c+ \trans{\AE}y + \trans{\AI}z),
\end{equation}
so we have
\begin{equation}\label{eq-theta}
  \theta(y,z)=\trans{\bE}\!y + \trans{\bI}\!z +\frac{1}{2}(\norm{c}^2-\norm{x(y,z)}^2).
\end{equation}

It is also not difficult to show \cite[Th.3.2]{malick-2004} that the
concave function $\theta$ is differentiable on $\RR^m$, and that its
gradient
\begin{equation}\label{eq-grad-theta}
  \nabla \theta(y,z)=-A x(y,z) + b
\end{equation}
is Lipschitz continuous. As any function with Lipschitz
gradient, $\theta$ is twice differentiable almost everywhere, but not
everywhere (this basically relies on the differentiability properties
of $\pK$; for the case $\Kk=\sdp$, see more in \cite{sun-sun-2002} and
\cite{malick-sendov-2005} among others).

The dual problem is thus
\begin{equation}\label{eq:dual}
\accol{\max\quad \theta(y,z)\\
\quad (y,z)\in \RR^{\mE}\times\RR_+^{\mI}.}
\end{equation}
Strong duality (the optimal values of \eqref{eq:pb2} and
\eqref{eq:dual} coincide) holds under a standard assumption in convex
optimization. The so-called (weak) Slater assumption (see
e.g.\,\cite{bertsekas-1995}, \cite{hull-1993}) is in our context:
\begin{equation}\label{eq:slater}
\exists \ \bar x \in \Pp \cap \Int \Kk.
\end{equation}
In fact, this assumption yields moreover
that there exists solutions to \eqref{eq:dual} (note that the
assumption has also a natural geometrical appeal in context of
projection methods, see \cite[Sec.\,3]{henrion-malick-2009}). Finally we
get directly the projection from dual solutions: let $(y^*,z^*)$ be a
(dual) solution of \eqref{eq:dual}, the (primal) solution $x^*$ of
\eqref{eq:pb} is the associated $x^*=x(y^*,z^*)$ (see \cite[Th.\,4.1]{malick-2004}).

\subsubsection{Apply standard algorithms}

To compute the projection of $c$ onto $\Pp\cap \Kk$, we just have to
solve the dual problem \eqref{eq:dual}. Let us have a closer look to
this problem: the constraints are simple positivity constraints on the
variable corresponding to the dualization of inequality constraints;
the dual function is a differentiable concave function with Lipschitz
gradient. This regularity has a huge impact in practice: it opens the
way for using standard algorithms for nonlinear optimization. Hence we
can use any of the following numerical methods to solve
\eqref{eq:dual} (as soon as the software can deal with the constraints
$z_i\geq 0$):
\begin{enumerate}
\item gradient methods: standard methods \cite{bertsekas-1995} or more
  evolved ones, as e.g.\;Nesterov's method \cite{nesterov-2004};
\item Newton-like methods: quasi-Newton, limited memory quasi-Newton,
  inexact Newton, Newton-CG, see textbooks \cite{nocedal-wright-1999}
  and \cite{bgls-2003} -- with the restriction that $\theta$ is not
  twice differentiable everywhere, so that we have to use the so-called
  semismooth Newton methods, see \cite{qi-sun-1993}.
\end{enumerate}
For example, \cite{malick-2004} uses a quasi-Newton method for solving
\eqref{eq:sdls}, and \cite{qi-sun-2006} uses a semismooth inexact
Newton method for solving \eqref{eq:corr}. We come back on these two
methods in the next section to give more practical details.

We also mention here the so-called inexact smoothing method of
\cite{gao-sun-2009} which consists in writing the optimality
conditions of the dual problem \eqref{eq:dual} as a nonsmooth fixed
point problem (and solving it by combining smoothing techniques and
an inexact Newton method; see e.g.\;\cite{nocedal-wright-1999}).

The dual problem \eqref{eq:dual} can thus be attacked with classical
tools or more evolved techniques. In practice, the choice of the
solving method depends on the structure of the problem and the target level
of sophistication. 

We call dual projection methods any method using an optimization code
for functions with Lipschitz gradient to maximize $\theta$ on
$\RR^\mE\times\RR_+^\mI$. Specifically, a dual projection method
generates a maximizing dual sequence $\{y_k,z_k\}_k$ together with the
primal sequence $x_k=x(y_k,z_k)$ such that:
\begin{eqnarray}
\theta(y_k,z_k) 
&=& \trans{\bE\!}\!y_k+\trans{\bI\!}\!z_k + \frac{1}{2}(\norm{c}^2-\norm{x_k}^2)\label{eq:sim}\\
\nabla \theta (y_k,z_k) &=& -Ax_k + b.\label{eq:simgrad}
\end{eqnarray}

We notice that in our numerical experiments with dual methods, we have
observed better behaviour and convergence when the (strong) Slater
assumption holds (that is, when \eqref{eq:slater} holds and moreover
$\AE$ is full rank).


\subsubsection{More algorithmic details (for the case without inequalities)}

We detail now further some algorithmic issues. To
simplify we focus on the case without inequalities ($\mI=0$, no dual
variables $z$). Iterations of most algorithms for maximizing $\theta$ can
be written as
\begin{equation}\label{eq:iter}
y_{k+1} = y_k + \tau_k W_k \nabla \theta (y_k).
\end{equation}
Note that the usual stopping test of these methods has an intrinsic
meaning: a threshold condition on the gradient
\begin{equation}\label{eq:stop-proj}
\norm{\nabla \theta(y_k)} = \norm{Ax_k-b} \leq \eps
\end{equation}
controls in fact the primal infeasibility.
Among these methods, let us discuss further the three following ones.

\paragraph{Gradient descent with constant step-size.}  We have a
  remarkable result: the gradient method in an adapted metric, namely
  \eqref{eq:iter} with 
  \begin{equation}\label{eq:grad}
    W_k = [A\trans{A}]^{-1}\qqandqq \tau_k =1,
  \end{equation}
  corresponds exactly to the alternating projection
  method~\eqref{eq:dykstra} (see \cite{malick-2004} for a proof in the
  special case of correlation matrices, and \cite{henrion-malick-2009}
  for the proof in general). We thus have a (surprizing) dual
  interpretation of the primal projection method. Using descent
  schemes more evolved than a simple gradient descent (see below) then
  leads to (dual) projection methods that can be seen as improvements
  of the basic alternating projection method.

\paragraph{BFGS Quasi-Newton method.} The method is known to be very
  efficient in general, and have many industrial applications (one of
  the most striking is in weather forecasting
  \cite{gilbert-lemarechal-1989}). The method can be readily applied
  to the dual problem, since it requires no more information than
  \eqref{eq:sim}: $W_k$ is constructed with successive gradients with
  the BFGS formula and $\tau_k$ is well-chosen with a Wolfe line-search (see
  e.g.\,\cite{bgls-2003}). The initial paper about dual methods
  \cite{malick-2004} proposes to use this method in general and
  reports very good numerical results on the nearest correlation
  matrix problem \eqref{eq:corr}. Since then, this dual method has
  been used successfully to solve real-life projection problems in
  numerical finance (among them: the problem of calibrating covariance
  matrices in robust portfolio selection \cite[5.4]{malick-2004}).  A
  simple Matlab implementation has been made publicly available
  together with \cite{henrion-malick-2009} for pedagogical and
  diffusion purposes.

\paragraph{Generalized (or semismooth) Newton.} A pure Newton method
  would be to use $\tau_k=1$ and $W_k = [H_k]^{-1}$ with the Hessian
  $H_k$ of $\theta$ at the current iterate $y_k$. In practice, an
  inexact generalized Newton method is used for the following
  reasons. 

  As mentioned earlier, $\theta$ is differentiable but not twice
  differentiable (though its gradient is Lipschitz continuous).  We
  can still replace the usual Hessian by a matrix $H_k\in
  \partial^2_c\theta(y_k)$ the Clarke generalized Hessian
  of $\theta$ at $y_k$~\cite{clarke-1983}.  Computing a
  matrix in $\partial^2_c\theta(y_k) \subset \sdp$ amounts to computing an element
  of the Clarke generalized Jacobian of the projection onto the cone
  $\partial_c\!\pK$ since we have (see \cite{hiriart-strodiot-nguyen-1984})
  \[
  \partial^2_c\theta(y_k) = A\,\partial_c\!\pK(c+\trans{A}\!y_k)\trans{A}.
  \]
  We can often compute an element of $\partial_c\!\pK$. For
  example, we even have an explicit expression of the whole
  $\partial_c\psdp$ \cite{malick-sendov-2005}.

  For overall efficiency of the method, the Newton direction $d_k$ is
  computed by solving the system $H_kd = \nabla\theta(y_k)$
  approximately, usually by conjugate gradient (CG) type
  methods. More precisely, the idea of so-called Newton-CG (also
  called inexact Newton going back to
  \cite{dembo-eisentat-steinhaug-1982}) is to stop the inner iteration
  of CG when
  \begin{equation}\label{eq-inexact}
    \norm{H_k d + \nabla \theta (\la_k)} 
    \leq \eta_k\norm{\nabla\theta(\la_k)}
  \end{equation}
  with small $\eta_k$ (see e.g.\,\cite{nocedal-wright-1999}). Note
  that preconditioning the Newton system is then crucial for
  practical efficiency. The nice feature of this algorithm is that
  $H_k$ has just to be known through products $H_kd$ so that
  large-scale problems can be handled. In our context, the main work on
  this method is \cite{qi-sun-2006} about the nearest correlation
  matrix problem; we come back to it in the next~section.

  We finish here with a couple of words about convergence of this
  Newton dual method. In general (see \cite{qi-sun-1993}), the two
  conditions to prove local superlinear convergence are that the
  minimum is strong (i.e.\;all elements of the generalized Hessian are
  positive definite), and the function has some smoothness (namely,
  the so-called semismoothness). In our situation, the two ingredients
  implying those conditions are the following ones:
  \begin{itemize}
    \item The intersection has some ``nondegeneracy'', in the sense of
      \cite[4.172]{bonnans-shapiro-2000} and
      \cite[Def.\,5]{alizadeh-haeberly-overton-1997}. This allows us to
      prove $\partial^2_c\theta(y_k) \succ 0$ (see
      e.g.\,\cite{qi-sun-2006} for a special case).
    \item The convex cone $\Kk$ has some ``regularity''. An example of
      sufficient regularity is that $\Kk$ is a semialgebraic set
      (i.e. defined by a finite number of polynomial (in)equalities).
      Indeed for semialgebraic convex sets, the projection $\pK$ and
      then $\theta$ are automatically semismooth
      \cite{bolte-daniilidis-lewis-2007} (which is the property needed to
      apply the convergence results of \cite{qi-sun-1993}. This is the case
      for direct products of the cones $\Ll_n$ and $\sdp$ (for which
      we even have strong semismoothness \cite{sun-sun-2002} so in
      fact quadratic convergence).
      \end{itemize}

\subsubsection{Illustration on nearest correlation matrix problem}

We give a rough idea of the efficiency of the dual approach on the
projection problem \eqref{eq:corr}. The first numerical results of
\cite[Sec.\,4]{malick-2004} show that the dual approach copes with
large-scale problems, in reporting that one can solve in a couple of
minutes projection problems of size around one thousand. By using the
dual generalized Newton method (instead of quasi-Newton as in
\cite{malick-2004}), the algorithm of \cite{qi-sun-2006}, improved
later by \cite{borsdorf-higham-2008}, gives very nice results in both
practice and theory. Nondegeneracy of the constraints and then of
the generalized Hessian is proved in \cite[Prop.\,3.6]{qi-sun-2006}: as
recalled above, this technical point leads to quadratic
convergence of the method \cite[Prop.\,5.3]{qi-sun-2006}.

Today's state of the art is that one can solve nearest correlation
matrix problems of big size (say, up to 4000-5000) in a reasonable
amount of computing time (say, less than 10 minutes on a standard
personal computer). The only limitation seems to be the memory
constraint to store and deal with dense large-scale data.

To give a more precise idea, let us report a couple of results from
\cite{borsdorf-higham-2008}.  The implementation of their dual
algorithm is in Matlab with some external Fortran subroutines (for
eigenvalues decomposition in particular). The stopping criterion is set to
\begin{equation}\label{eq:approx}
  \norm{\nabla \theta(y_k)} \leq 10^{-7}n.
\end{equation}
We consider the nearest correlation matrix problems for two (non-SDP)
matrices with unit diagonal (of size $n_1=1399$ and $n_2= 3120$)
provided by a fund management company. The dual method solves them in
around 2 and 15 min., respectively, on a very standard machine (see
more details in \cite{borsdorf-higham-2008}).

We finish with a last remark about accuracy. The approximate
correlation matrix $X$ that is computed by such a dual method is often
just what is needed in practice. It might happen though that a special
application requires a perfect correlation matrix -- that is, with exactly
ones on the diagonal, whereas $X$ satisfies only (by \eqref{eq:approx})
\[
\Big(\sum_{i=1}^n(X_{ii}-1)^2\Big)^{-1/2}\leq 10^{-7}n.
\]
A simple post-treatment corrects this. Setting diagonal
elements to ones may destroys the positiveness, so we apply the usual
transformation that computes the associated correlation matrix $\bar X$
from a covariance matrix $X$, namely
\[
\bar X = D^{-1/2}XD^{-1/2} \qqandqq D=\diag(X).
\]
This operation increases the distance from $C$; but the error is still
under control (by $\eps/(1-\eps)$; see
\cite[Prop.\;3.2]{borsdorf-higham-2008}).

\subsection{Discussion: applications, generalizations}

\subsubsection{Direct or indirect applications}

Conic projection problems with the positive semidefinite cone (like
$\Kk=\sdp$, $\Kk=\sdp\times (\RR^+)^p$ or $\Kk=\Ss^+_{n_1}\!\times
\cdots \times \Ss^+_{n_p}$) are numerous in engineering. Constructing
structured semidefinite matrices, for example, are naturally modeled
this way. Such problems naturally appear in finance for constructing
structured covariance matrices (as a calibration step before
simulations); they also appear in many other fields, such as in control
(e.g.\,\cite{ieee-2009}), in numerical algebra
(e.g.\,\cite{hankel-2007}), or in optics
(e.g.\,\cite{vandenberghe-2008}), to name a few of them.

Conic projections also appear as inner subproblems within more
involved optimization problems. Solving efficiently these inner
problems is often the key to numerical efficiency of the overall
approach. Let us give some examples.
\begin{itemize}
\item {\em Linear conic programming.} So-called regularization methods
  for solving~\eqref{eq:lin} use the conic projection problem as an inner
  subproblem; these methods are studied in \secref{sec:reg}.

\smallskip

\item {\em Weighted projections.} For given weights $H_{ij}\geq
  0$, consider the semidefinite projection~\eqref{eq:sdls} with a
  different objective function
  \[
    \accol{\min \quad \frac{1}{2}\sum^n_{i,j=1}H_{ij}(X_{ij}-C_{ij})^2\\ 
      \quad \prods{A_i}{X}= (\textor \leq)\;b_i, \ \ i=1,\ldots,m\\ 
      \quad X\succeq 0.}
  \]
  An augmented Lagrangian approach for this problem \cite{qi-sun-2010}
  produces a projection-like inner problem, which is solved by
  a semismooth Newton method (recall the discussion of the previous
  section).

\smallskip
\item {\em Low-rank projections.}  Consider the semidefinite
  projection problem~\eqref{eq:sdls} with additional rank-constraint
  \begin{equation}\label{eq:low}
    \accol{\min \quad \frac{1}{2}\norm{X-C}^2\\ 
      \quad \prods{A_i}{X}= (\textor \leq)\;b_i, \ \ i=1,\ldots,m\\ 
      \quad X\succeq 0, \ \rank X =r.}
  \end{equation} 
  This difficult non-convex calibration problem has several
  applications in finance and insurance industry (e.g.\;pricing
  interest rate derivatives for some models; see
  e.g.\,\cite{brigo-mercurio-2006}). Two approaches (by augmented
  Lagrangian \cite{li-qi-2010} and by penalty techniques
  \cite{gao-sun-2010}) have been recently proposed to solve these types
  of problems; both approaches solve a sequence of projection-like
  subproblems. The numerical engine is a dual semismooth truncated
  Newton algorithm for computing projections.
\end{itemize}

For these applications of conic projections, the techniques and the
arguments are often the same, but are redeveloped for each particular
projection problem encountered.  We hope that the unified
view of Section 2 can bring forth the common ground of these
methods and to better understand how and why they work well. We
finish this section by pointing out an easy geometrical application.

\subsubsection{Application for solving conic feasibility problems}

The conic feasibility problem consists simply in finding a point $x$
in the intersection $\Kk\cap\Pp$. Many engineering problems can be
formulated as semidefinite or conic feasibility problems (for example
in robust control \cite{boyd-1994} where an element in the
intersection is a certificate of stability of solutions of
differential equations).  \secref{sec:SOSfeas} focuses on semidefinite
feasibility problems arising when testing positivity of
polynomials. We refer to the introduction of
\cite{henrion-malick-2009} for more examples and references.

A simple and natural technique for solving conic feasibility problems
is just to project a (well-chosen) point onto the intersection
$\Kk\cap\Pp$ (by dual projection methods for example).  In
\cite{henrion-malick-2009}, a comparative study of such a conic
projection method with the usual approach using SeDuMi was carried out
precisely on polynomial problems. It was shown there that an
elementary Matlab implementation can be competitive with a
sophisticated primal-dual interior-point implementation. This would
even have a better performance if an initial heuristic for finding a
good point to project could be determined (the numerical experiments
of \cite[Sec.\,6]{henrion-malick-2009} simply use $c=0$). An answer to
this latter point is provided by the regularization methods of the
next section.

\section{Projections in regularization methods}\label{sec:reg}

We focus in this section on standard linear conic programming. We show
that, following classical convex optimization techniques, conic
projections can be used to solve linear conic programming problems.

There exist many numerical methods for solving linear conic problem
\eqref{eq:lin} (see the first handbook
\cite{wolkowicz-saigal-vandenberghe-2000}).  But on the other hand,
there also exist big conic problems, and especially big SDP problems,
that make all the standard methods fail. Relaxations of combinatorial
optimization problems and polynomial optimization problems yield
indeed challenging problems. This motivates the development of new
algorithmic schemes.

The strategy that we present in this section exploits the
efficiency of projection methods by developing proximal algorithms for
linear conic programming.  We generalize the developments of
\cite{malick-povh-rendl-wiegele-2009}, and give all way long
references to related works. As for numerical aspects, the target
problems are semidefinite programs with the number of constraints
possibly very large (more than $100,\!000$).

\subsection{Proximal method for linear conic programming}\label{sec:prox}

\subsubsection{Apply classical techniques of convex optimization}

The proximal algorithm is a classical method of convex optimization
and variational analysis: it goes back from the 1970s with premises in
\cite{bellman-kalaba-lockett-1966}, the first work
\cite{martinet-1970} and the important reference
\cite{rockafellar-1976}.  The driving idea of the proximal algorithm
is to add quadratic terms to the objective function to ``regularize'' the
problem (ensuring existence, uniqueness, and stability of solutions).
A (primal) proximal method of the linear conic problem \eqref{eq:lin}
goes along the following lines. 

Consider the problem with respect to $(x,p)$
\[
\accol{\min \quad \trans{c}{x} + \frac{1}{2t}\norm{x-p}^2\\
  \quad p\in \RR^n, \ x\in \Pp\cap \Kk.}
\]
By minimizing first with respect to $p$, we see that this problem is
equivalent to the primal linear conic problem \eqref{eq:lin}. We have
added to the objective function a quadratic ``regularizing'' term
$\norm{x-p}^2$ with the so-called ``prox-parameter''~$t$.  The idea
now is to solve this problem in two steps: first with respect to $x$,
and second to $p$:
\begin{equation}\label{eq:proxform}
\begin{array}{rl}
\accol{
  \min\\{\quad p \in \RR^n}\quad}
&
\!\!\!\!\!\!\!\!\left(\begin{array}{l}\min\quad\trans{c}{x} + \frac{1}{2t} \norm{x-p}^2\\
  \quad x\in \Pp\cap\Kk\end{array}\right).
\end{array}
\end{equation}
The outer problem is thus the minimization with respect to $p$
of the function 
\begin{equation}\label{eq:F}
F(p):=\accol{\min \quad \trans{c}{x} + \frac{1}{2t}\norm{x-p}^2\\
  \quad x\in \Pp\cap \Kk}
\end{equation}
which is the result of the inner optimization problem parametrized
by~$p$.  As such defined, $F$ is the so-called Moreau-Yosida
regularization of the function $x\ra \trans{c}x + i_{\Pp\cap\Kk}(x)$ the
linear objective function plus the indicator function of the
intersection (see e.g.\,\cite[Ch.XV\!.4]{hull-1993}).

The connection with the previous developments of this chapter is then
obvious: the above inner problem is essentially a projection problem
as studied in the previous section (see \eqref{eq:gene}). The solution
of the inner problem (the ``projection'') is called the proximal point
and denoted
\[
\prox(p):=\left\{\begin{array}{l}\argmin\quad\trans{c}{x} +
\frac{1}{2t} \norm{x-p}^2\\ \quad x\in \Pp\cap\Kk.\end{array}\right.
\]
Note that, for simplicity, the dependence of $F$ and $\prox$ with
respect to $t$ is dropped in notation.

\subsubsection{Primal proximal algorithm for conic programming}

Applying basic convex analysis properties, it is easy to prove (see
e.g.\,\cite[Ch.XV\!.4]{hull-1993}) that the Moreau-Yosida
regularization $F$ is convex and differentiable with gradient $\nabla
F(p) = (p-\prox(p))/t$. The optimality condition of the unconstrained
minimization of $F$ is then simply
\begin{equation}\label{eq:fix}
\bar p=\prox(\bar p).
\end{equation}
Moreover a fixed-point $\bar p$ of the $\prox$ operator is
also a solution of the initial linear conic problem~\eqref{eq:lin}:
observe indeed that $\bar p$ is feasible and reaches the optimal
value, since
\begin{equation}\label{eq:sol}
\val\eqref{eq:lin} = \min F(p) = F(\bar p) = \trans{c}\bar p.
\end{equation}

A (primal) proximal algorithm for solving \eqref{eq:lin} then consists
of a fixed-point algorithm on \eqref{eq:fix}
\begin{equation}\label{eq:proxiter}
p_{k+1}=\prox(p_k).
\end{equation}
Since computing $\prox(p_k)$ corresponds to solving a projection
problem, we can use any of the algorithmic schemes described in
\secref{sec:algo} to implement~\eqref{eq:proxiter} inside of this proximal
algorithm.  We call the proximal method the outer algorithm, and the
chosen projection algorithm, the inner algorithm. We study in
\secref{sec:reg} the family of proximal algorithms obtained when dual
projection algorithms of \secref{sec:dual} are used as inner
algorithms.

As we have an iterative optimization algorithm (inner algorithm)
inside of another iterative algorithm (outer algorithm), the question
of the stopping tests of the inner algorithm is obviously crucial. For
practical efficiency, the inner stopping test should somehow depend on
some outer information; we~come back later in detail to this important
point.

So in fact the iteration \eqref{eq:proxiter} is not carried out
exactly, and replaced instead by a looser implementable relation
\begin{equation}\label{eq:proxiterimpl}
\norm{p_{k+1}-\prox(p_k)}\leq \eps_k.
\end{equation}
Whatever is the inner projection algorithm, we have the general
global convergence of the method under the assumption that the
sequence of errors $\eps_k$ goes rapidly to zero.

\begin{proposition}[Global convergence]\label{prop:convergence}
  Assume that there exist a solution~to~\eqref{eq:lin}.  If $(t_k)_k$
  is bounded away from $0$ and if the primal proximal algorithm generates a
  sequence $(p_k)_k$ such that
  \begin{equation}\label{eq:sum}
  \sum_k\eps_k < + \infty
  \end{equation}
  then $(p_k)_k$ converges to a solution $\bar p$ of \eqref{eq:lin}.
\end{proposition}

{\bf Proof:}
  The result is straightforward from the general convergence result of
  proximal algorithms. As a consequence of \eqref{eq:sum} and the
  existence of a solution to \eqref{eq:lin}, the sequence $(p_k)_k$ is
  bounded and we can apply \cite[Th.1]{rockafellar-1976}: $(p_k)_k$
  converges to a fixed-point to $\prox$ which is a solution of
  \eqref{eq:lin} by \eqref{eq:sol}. $\Box$

\subsubsection{Dual point of view: augmented Lagrangian}

We give here some details about the dual interpretation of the above
primal algorithmic approach. It is known indeed that a proximal method
for a problem corresponds exactly to an augmented Lagrangian method on
its dual; we~detail this for our case. To simplify writing duals, we
abandon the general formulation \eqref{eq:lin}, and we suppose that
there is no affine inequalities (or that there are incorporated with
slack variables in $\Kk$). So we work from now with the standard form
of primal and dual linear conic problems
\begin{equation}\label{eq:primaldual}
\accol{\min \quad \trans{c}x\\
\quad Ax=b \\
\quad x\in \Kk}
\qqandqq
\accol{\max \quad \trans{b}y\\
\quad \trans{A}y-u-c =0\\
\quad u\in \Kk^o.}
\end{equation}

Augmented Lagrangian methods are important classical regularization
techniques in convex optimization (see \cite{polyak-tret-1972},
\cite{rockafellar-1976b} for important earlier references, and
\cite[Chap.XII]{hull-1993} for the connection with usual Lagrangian
duality). In our situation, a dual augmented Lagrangian method goes along
the following lines. Introduce the augmented Lagrangian function $L$
with parameter $t >0$, for the dual problem \eqref{eq:primaldual}:
$$
L(y,u;p) := \trans{b}y - \trans{p}(\trans{A} y - u -c) -
\frac{t}{2}\norm{\trans{A} y - u -c}^2.
$$ 
Note that this is just the usual Lagrangian for the problem
\begin{equation}\label{eq-dual-stab}
  \accol{
    \max\quad \trans{b}y - \frac{t}{2}\norm{\trans{A} y - u -c}^2\\
    \quad \trans{A}y - u -c = 0, ~ u \in \Kk^o,
  }
\end{equation}
that is the dual problem with an additional redundant quadratic
term in the objective. The convex (bi)dual function is then defined as
\begin{equation}\label{eq-inner-rendl}
  \Theta(p):= \max_{y\in \RR^m,u\in \Kk^o} L(y,u;p).
\end{equation}
The bridge between the primal proximal method and the dual augmented
Lagrangian is set in the next proposition, formalizing
a well-known result.
\begin{proposition}
  With notation above, we have $\Theta(p)=F(p)$ for $p\in \RR^n$.
\end{proposition}

{\bf Proof:}
  Just apply \cite[XII.5.2.3]{hull-1993}: the augmented Lagrangian function 
  $\Theta(p)$ is the Moreau-Yosida of the usual dual function, 
  which is here
  \[
  \trans{c}p + i_{\{Ax=b\}\cap\Kk}(p) = 
  \max_{y,u\in \Kk^o}\trans{b}y - \trans{p}(\trans{A} y - u -c).
  \]
  This is exactly $F(p)$ defined by \eqref{eq:F} (in the case when
  $\Pp$ is just the affine subspace of equation $Ax=b$).$\Box$

The primal regularization by proximal approach and the dual augmented
Lagrangian regularization thus correspond exactly to the same
quadratic regularization process viewed either on the primal problem
or on the dual~\eqref{eq:primaldual}.

The developments of this section share similar properties with
other augmented Lagrangian-type approaches for conic
programming, among them: a primal augmented Lagrangian in
\cite{burer-vandenbusshe-2006}, a primal-dual augmented Lagrangian in
\cite{jarre-rendl-2007} and a penalized augmented Lagrangian in
\cite{kocvara-stingl-2007}.

\subsection{Regularization methods for linear conic programming}

In this section we give more details on primal proximal algorithms (or
dual augmented Lagrangian algorithms) that use dual projection methods as
inner algorithms to carry out the proximal iteration
\eqref{eq:proxiterimpl}.  This family of algorithms is introduced for
the case $\Kk=\sdp$ in \cite{malick-povh-rendl-wiegele-2009}. They are
called regularization algorithms (rather than proximal algorithms,
which would focus on the primal point of view only); we keep this terminology
here. This section is more technical and could be skipped at a first
reading.

Regularization algorithms for conic programming specialize on
three points:
\begin{enumerate}
\item the dual projection algorithm to compute $\prox(x_k)$,
\item the rule to stop this inner algorithm,
\item the rule to update the prox-parameter $t_k$.
\end{enumerate}
The third point is an inherent difficulty of any practical
implementation of proximal methods (e.g.\;bundle methods, see
\cite{correa-lemarechal-1993}). We are not aware of general techniques
to tackle it. So we focus here on the first two points.

\subsubsection{Dual projection methods as inner algorithms}

We could use any dual projection algorithm of \secref{sec:dual} to solve
\begin{equation}\label{eq:inner}
\accol{\min \quad \trans{c}{x} + \frac{1}{2t}\norm{x-p}^2\\
  \quad Ax= b, \ x\in \Kk.}
\end{equation}
Embedded in a proximal scheme, a dual projection algorithm would lead
to the forthcoming overall algorithm for solving linear conic problems
\eqref{eq:primaldual}.

Note first that equations \eqref{eq-xy} and \eqref{eq-theta}
for the projection-like problem \eqref{eq:inner}
become respectively
\begin{eqnarray}\label{eq:xy2}
x(y)         &=& \pK\big(p+t(\trans{A}y-c)\big)\\
\theta(y)    &=& \trans{b}\!y+\frac{1}{2t}(\norm{p}^2-\norm{x(y)}^2).
\end{eqnarray}
We use the (slightly loose) formulation \eqref{eq:iter} of the
iteration of dual projection methods to write a general
regularization algorithm. We index the outer iterations by $k$ and the
inner ones by $\ell$.

\begin{algorithm}[Regularization methods]\hfill
  
  Outer loop on $k$ stopped when  $\norm{p_{k+1}- p_k}$ small:
   
   \quad Inner loop on $\ell$ stopped when $\norm{Ax_\ell-b}$ small enough:

   \quad \quad Compute $x_\ell = \pK(p_k + t_k(\trans{A} y_\ell-c))$ and  $g_\ell = b - Ax_\ell$
  
   \quad \quad Update $y_{\ell+1} =  y_\ell + \tau_\ell\, W_\ell\, g_\ell$ with appropriate $\tau_\ell$ and $W_\ell$
  
   \quad end (inner loop) 
  
   \quad Update $p_{k+1} = x_\ell$ (and $t_k$)

   end (outer loop)
   \end{algorithm}
  
\noindent We discuss several points about the above conceptual
algorithm.

\begin{itemize}

\item \emph{Memory.}  An important feature of regularization methods is
  the rather low memory requirement. The intrinsic operations of the
  algorithm are basically the projection onto the cone and the
  multiplications by $A$ and $\trans{A}$. If the data has some
  structure, those multiplications could be performed efficiently
  (without constructing matrices). Moreover for maximizing
  $\theta$ (that is, essentially, implementing $y_{\ell+1} = y_\ell +
  \tau_\ell\, W_\ell\, g_\ell$), we could use algorithms of smooth
  unconstrained optimization adapted to large-scale problems and then
  requiring low-memory (as limited memory BGFS or Newton-CG, see
  e.g.\;\cite{nocedal-wright-1999} and \cite{bgls-2003}). We come back
  to this point later when discussing numerical issues. Roughly
  speaking, the point is the following: the computer memory can
  be used for storing problem data and the computation does not
  require much more extra memory.

\smallskip

\item \emph{Inner restarting.} At the outer iteration $k$, the inner
  projection algorithm can be initialized with the best $y_\ell$ of
  the previous iteration $k-1$. This has an intuitive appeal, so that
  in practice, $\ell$ keeps increasing over the outer
  iterations. (Note also that the historical information on gradients may
  be carried out from iteration $k-1$ to $k$ as well.)

\smallskip

\item \emph{Dual variable $u$.} It is known that for any $x\in \RR^n$,
  the projection onto the polar cone $\proj_{\Kk^o}(x)$ is
  given by $\pK(x)+\proj_{\Kk^o}(x)=x$ (together with
  $\trans{\pK(x)}\proj_{\Kk^o}(x)=0$, see
  \cite[III.3.2.5]{hull-1993}).  When computing $x_\ell$, we thus get
  automatically
  \[
  u_{\ell} = \proj_{\Kk^o}(p_k + t_k(\trans{A} y_\ell-c))/t_k
  \]
  and it holds
  \begin{equation}\label{eq:ul}
  p_k + t_k(\trans{A} y_\ell-c) = t_ku_\ell + x_\ell.
  \end{equation}

\smallskip  

\item \emph{Dual outer iterates.}  At the end of outer iteration $k$,
  we set (with a slight abuse of notation) $y_{k+1}=y_\ell$ and
  $u_{k+1}=y_\ell$ for $\ell$ the final iteration of inner
  algorithm. Thus we have a sequence of primal-dual outer iterates
  $(p_k,y_k,u_k)\in \Kk\times\RR^n\times\Kk^o$.  Under some technical
  assumptions, we can prove a convergence result of the same vein as
  \propref{prop:convergence}: any accumulation point of the sequence
  $(p_k,y_k,u_k)$ is a primal-dual solution of \eqref{eq:pb2} (see
  e.g.\;Theorem\,4.5 of \cite{malick-povh-rendl-wiegele-2009} for a
  proof when $\Kk=\sdp$).

\smallskip

\item \emph{Outer stopping test.} We have already noticed in
  \eqref{eq:stop-proj} that the natural stopping test of dual
  projection algorithms controls primal infeasibility
  $\norm{Ax_\ell-b}$.  Interpreted as a fixed point iteration
  \eqref{eq:proxiter}, the natural stopping of the proximal algorithm
  is $\norm{p_{k+1}-p_k}$; it turns out that this can interpreted as
  controlling dual infeasibility.  Note indeed that \eqref{eq:ul}
  yields
  \[
  p_k + t_k(\trans{A} y_k-c) = t_ku_k + p_k
  \]
  and then we have
  \[
  \norm{p_{k+1}-p_k} = t_k\norm{\trans{A} y_k-u_k-c}.
  \]

\smallskip

\item \emph{Normal to interior-point methods.} By construction, conic
  feasibility $p\in \Kk$ (and $u\in \Kk^o$) and complementary
  $\trans{x}u=0$ are ensured throughout the algorithm, while
  primal-dual feasibilities are obtained asymptotically. In~contrast,
  recall that basic interior-point methods maintain primal and dual
  feasibility and the conic feasibility and work to reach
  complementarity. Note also that, regularization algorithms give
  solutions that are usually on the boundary of the cone $\Kk$, since
  the primal iterates are constructed by projections onto $\Kk$. In
  contrast again, basic interior-points give solutions as inside of
  the cone as possible. In a sense, regularization methods are then
  ``normal'' to interior point methods.

\smallskip

\item \emph{Overall stopping test.} We have seen above that the natural
  outer and inner stopping rules of the regularization algorithm have
  a practical interpretation as dual and primal infeasibilities. Since
  complementary and conic feasibility are ensured by construction, the
  natural stopping test of the overall algorithm is
  \begin{equation}\label{eq:residual}
  \max\left\{\norm{Ap_k-b}, \norm{\trans{A} y_k-u_k-c}\right\}.
  \end{equation}
  In practice, one should divide moreover the two infeasibilities by
  some constant quantities to get homogeneous ratios.
\end{itemize}

\subsubsection{Stopping inner iterations}

For which inner iteration $\ell$ can we set $p_{k+1}=x_\ell$ to
proceed with the outer iteration? Since we have a loop inside of
another, the rule to terminate the inner projection algorithm is
indeed an important technical point for regularization algorithms.  We
discuss three strategies to set up inner stopping~rules.

\paragraph{Solving approximately the inner problem.}  

The usual
  stopping inner rule in proximal methods is to stop inner iterations
  when the current inner iterate $x_\ell$ is close to the proximal
  point $\prox(p_k)$. Doing this, the regularization algorithm
  approximates at best the conceptual proximal algorithm (which
  requires to solve the inner problem exactly), so that we keep
  convergence properties (as in \propref{prop:convergence} for
  instance).

  This is the strategy followed by the regularization method of
  \cite{zhao-sun-toh-2010} for $\Kk=\sdp$. This paper adopts the dual
  point of view (augmented Lagrangian) and uses semismooth Newton as
  dual projection algorithm for inner iterations (remember
  \secref{sec:dual}). The regularization method thus combines the
  usual stopping strategy and an efficient inner algorithm (supporting
  large-scale problems); it gives excellent numerical results on
  various semidefinite programming test-problems (see the last section
  of \cite{zhao-sun-toh-2010}). Under nondegeneracy assumptions, we
  have moreover proofs of global and local convergences of the inner
  algorithm as well as the overall regularization method.

\paragraph{Only one inner iteration; interpretation as saddle-point.}

  An opposite strategy is to do only one inner iteration per outer
  iteration. This cautious strategy is motivated by the
  following remark. Let us come back to the proximal formulation
  \eqref{eq:proxform} of the linear conic problem.  Under Slater
  assumption \eqref{eq:slater}, the inner projection problem can be
  replaced by its dual (recall \eqref{eq-theta}, \eqref{eq:dual} and
  \eqref{eq:xy2}), so that the primal and dual conic problems
  \eqref{eq:primaldual} have the saddle-point formulation
  \[
  \begin{array}{rl}
    \accol{ \min\\{\quad p \in \RR^n}\quad} &
    \!\!\!\!\!\!\!\!\left(\begin{array}{l}\max \quad \trans{b}{y} - 
      \frac{1}{2t}(\norm{p}^2 - \norm{\pK(p + t(\trans{A} y -c)}^2) \\ \quad y\in
      \RR^m \end{array}\right).
  \end{array}
  \]
  With this point of view on the process, the choice of inner stopping
  conditions appears indeed to be crucial, because the inner and
  outer loops are antagonistic, as the first minimizes and the second
  maximizes. The idea of the ``Arrow--Hurwicz'' approach (see, for
  instance, \cite{arrow-hurwicz-uzawa-1959}) is essentially to proceed
  with gradient-like iterations with respect to each variable
  successively.

  This is the strategy of the simple regularization method presented
  in \cite[Sec.\,5.1]{malick-povh-rendl-wiegele-2009}. Doing one inner
  iteration of \eqref{eq:iter} with $W_k=[A\trans{A}]$ and
  $\tau_k=1/t_k$ allows to simplify the regularization algorithm to
  just one~loop~with
  \begin{eqnarray}
  p_{k+1} &=& \pK(p_k + t_k(\trans{A} y_\ell-c)) 
  \quad (\text{with $u_{k+1}$ as by-product})\label{eq:pk}\\
  y_{k+1} &=&  y_k + [A\trans{A}]^{-1}(b - Ap_k)/t_k.\label{eq:yk}
  \end{eqnarray}
  We retrieve algorithm 5.1 of \cite{malick-povh-rendl-wiegele-2009}
  by using the proof of proposition 3.4 in there.  This simple
  regularization algorithm has also an interpretation as an
  alternating direction method, see the chapter of this book devoted
  to them.

  In practice, it is important to note that $A\trans{A}$ and its
  Cholesky factorization can be computed only once at the
  beginning of the algorithm. Even if this is an expensive task in
  general for problems that have many unstructured constraints (so
  that $A\trans{A}$ is big and unstructured), there exists some cases
  when $A\trans{A}$ is very simple, so that solving the system
  \eqref{eq:yk} is cheap. This is the case in particular when
  $A\trans{A}$ is diagonal, as for SDP relaxations of max-cut problem,
  or $k$-max-cut problem \cite{goemans-williamson-1995}, frequency
  assignment problems, see \cite[(5)]{burer-monteiro-zhang-2003},
  max-stable problem, see more below, and polynomial minimization
  problems, see forthcoming \secref{sec:poly}.
  
\paragraph{Something in-between.} An attractive option is to find
  something in-between the previous two extreme strategies. Explicit
  rules for the management of $\eps_k$ in \eqref{eq:proxiterimpl}
  should be given, and for numerical efficiency they should be given
  online. Using off-line rules independent of the function is
  interesting in theory since it allows to get proof of linear
  convergence. Following usual paradigms of numerical optimization, it
  would be possible to do better as soon as practical efficiency in
  concerned. An appropriate stopping rule still has to be found and
  studied; this is actually a general question and we are not aware of
  efficient techniques.

  Note finally that the practical implementation of the regularization
  method of \cite{zhao-sun-toh-2010} does indeed something in-between:
  the simple scheme with one inner gradient iteration (second
  strategy) is used as a preprocessing phase before switching to
  making inner Newton iterations (first strategy). See more about this
  on numerical illustrations of the method in \secref{sec:poly}. A
  stopping test along the above lines would further enhance the
  numerical performance of the implementation.


\subsubsection{Illustration: Computing Lov\'asz theta number}

We finish this section with a numerical illustration (borrowed from
\cite{malick-povh-rendl-wiegele-2009}) of the performance of
regularization methods on a classical combinatorial optimization
problem.

Lov\'asz \cite{lovasz-1979} proved the celebrated ``sandwich'' theorem
in graph theory: the stable number $\alpha(G)$ of a graph $G$ and the
chromatic number $\chi(\bar G)$ of its complementary graph $\bar G$
are separated
\[
\al(G) \leq \vartheta(G) \leq \chi(\bar G)
\]
by the optimal value of an SDP problem
\begin{equation}\label{eq:lovasz}
\vartheta(G)= 
\accol{\max\quad \langle \mathbf{1}_{n\times n}, X \rangle \\
  \quad X_{ij}=0,\quad \text{when $(i,j)$ is an edge of $G$}\\
  \quad \trace\,X=1, \ \ X \succeq 0.}
\end{equation}
As expected, it can be shown that this SDP problem is a formulation of
the SDP relaxation of the max-stable problem.  The stable number
$\alpha(G)$ and the chromatic number $\chi(\bar G)$ are both NP-hard
to compute and even hard to approximate, so that the tractable
$\vartheta(G)$ gives interesting information about~$G$.
  
Some graphs from the DIMACS collection~\cite{johnson-trick-1996} are
very challenging instances for computing $\vartheta(G)$. Those graphs
have many edges and also many edges on the complementary, so that
makes them the most difficult for standard methods (as noticed in
\cite{dr:07}). On the other hand, the structure of problem
\eqref{eq:lovasz} is very favorable for regularization methods and in
particular for the simple one of
\cite[Sec.\,5.1]{malick-povh-rendl-wiegele-2009} which is essentially
\eqref{eq:pk}-\eqref{eq:yk}.  Observe indeed that the affine
constraints \eqref{eq:lovasz} is given by ``orthogonal'' matrices
(i.e.\;$\prods{A_j}{A_i}=0$), such that the matrix $A\trans{A}$ is
diagonal. For illustration, the next table reports some of the bounds
that were computed for the first time in
\cite{malick-povh-rendl-wiegele-2009}.
\begin{center}
  \begin{tabular}{l|c|c|c}
    graph name & $n$ & $m$ & $\vartheta(G)$ \\
    \hline
    brock400-1   & 400 &  59723 & 10.388 \\
    keller5      & 776 & 225990 & 31.000 \\
    brock800-1   & 800 & 207505 & 19.233 \\
    p-hat500-1   & 500 &  31569 & 58.036 \\
    p-hat1000-3  &1000 & 371746 & 18.23 \\
  \end{tabular}
  \end{center}
For more examples, see
\cite[Sec.\;5.4]{malick-povh-rendl-wiegele-2009} and
\cite[Sec.\;6.3]{zhao-sun-toh-2010}.

\section{Applications to polynomial optimization}\label{sec:poly}

In this section, we illustrate the regularization methods for solving
linear semidefinite optimization problems in the context of polynomial
optimization. We collect numerical experiments showing that
regularization algorithms can be considered as an alternative to
standard methods for deciding whether a polynomial is non-negative
(Section \ref{sec:SOSfeas}) and for globally minimizing a polynomial
(Section \ref{sec:SOSmin}). The point is thus the same as
in \cite{nie-2009} which reports extensive numerical results using
regularization methods for solving various large-scale polynomial
optimization problems.

We aim at giving here a methodology: our main focus is the generality
of the approach and the reproducibility of experiments and results. We
explain how to generate the test problems, we use public-domain
implementations of the algorithms with default parameter tunings and
with no attempt to adapt them to each particular problem instances
(contrary to \cite{nie-2009}). We do not carry out a comprehensive
benchmarking with all methods and solvers; we just compare a widely
used implementation of interior-point algorithm
\cite{sturm-1999} with two recent implementations of regularization
methods (the basic one of
\cite[Sec.\,5.1]{malick-povh-rendl-wiegele-2009}, and the more
sophisticated one of \cite{zhao-sun-toh-2010}).

\subsection{Sum-of-squares, SDP and software}\label{sec:SOSintro}

We briefly introduce in this section the notions and notation of
polynomial optimization that we need. We refer to the recent surveys
\cite{laurent-2009} and \cite{lasserre-2009} and to the other chapters
of this book for more on this topic.

Consider a multivariate polynomial of total degree $2d$
\begin{equation}\label{eq:poly}
v \in {\mathbb R}^N \longmapsto p(v) = \sum_{|\alpha|\leq 2d} p_{\alpha} v^{\alpha}.
\end{equation}
We use here the multi-index notation $v^{\alpha} = v^{\alpha_1}_1
\cdots v^{\alpha_N}_N$ where $\alpha \in {\mathbb N}^N$ runs over all
nonnegative integer vectors of sum $|\alpha| = \alpha_1 + \cdots +
\alpha_N \leq 2d$.  We say that $p(v)$ is a sum-of-squares (SOS) of
polynomials if one could find polynomials $q_k(v)$~such~that
\begin{equation}\label{eq:SOSpoly}
p(v) = \sum_k {q_k}^2(v).
\end{equation}
It can be shown that finding such polynomials $q_k(v)$ amounts to a
semidefinite feasibility problem. More specifically, if $\pi(v)$
denotes a vector of basis of polynomials of total degree less than or equal to
$d$, finding an SOS decomposition (\ref{eq:SOSpoly}) amounts to finding
a so-called Gram matrix $X \in {\mathbb R}^{n\times n}$ such that
\begin{equation}\label{eq:SOSGram}
p(v)=\trans{\pi(v)} X \pi(v) \qqandqq X\in \sdp.
\end{equation}
The set of SOS polynomials has thus a SDP representation of the form
\begin{equation}\label{eq:Gramvector}
Ax=b, \quad x \in {\mathcal K}
\end{equation}
where $\mathcal K =\sdp$, $A \in {\mathbb R}^{m\times n^2}$ is a
linear operator depending only on the choice of basis $\pi(v)$, and $b
\in {\mathbb R}^m$ is a vector depending on $p(v)$.  For example if
the vector of basis polynomials $\pi(v)=[x^{\alpha}]_{|\alpha|\leq d}$
contains monomials $x^{\alpha}$ indexed by $\alpha \in
{\mathbb N}^n$, then identifying powers of $v$ in relation
(\ref{eq:SOSGram}) yields
\[
p_{\alpha} = \prods{A_{\alpha}}{\pi(v)\trans{\pi(v)}}, \quad \text{ for all } \alpha
\]
where matrix $A_{\alpha}$ selects monomials $x^{\alpha}$ in rank-one
matrix $\pi(v)\trans{\pi(v)}$. More specifically, the entry in $A_{\alpha}$
with row index $\beta$ and column index $\gamma$ is equal to one if
$\beta+\gamma=\alpha$ and zero otherwise.  In problem
(\ref{eq:Gramvector}), the row indexed by $\alpha$ in matrix $A$
collects entries of matrix $A_{\alpha}$, and the row indexed by
$\alpha$ in vector $b$ is equal to $p_{\alpha}$. Note that
\[
n = \left(\begin{array}{c}N+d\\N\end{array}\right) \qqandqq
m = \left(\begin{array}{c}N+2d\\N\end{array}\right),
\]
so that the sizes of SDP problems grow quickly with the degree and the
number of variables. 

The important remark is that this type of constraints are favorable to
regularization methods: $A\trans{A}$ is always diagonal indeed. To see
this, let $\alpha$, $\beta$ denote the row and column indices in
matrix $A\trans{A}$. By construction, the entry $(\alpha,\beta)$ in
$A\trans{A}$ is equal to $\prods{A_{\alpha}}{A_{\beta}}$: if
$\alpha=\beta$, this is equal to the number of non-zero entries in
matrix $A_{\alpha}$, otherwise, this is zero.  Since it is important
for numerical efficiency, we formalize the previous remark in a
proposition.

\begin{proposition}[Orthogonality of constraints]
  Let $A$ be the matrix in SOS semidefinite problem
  (\ref{eq:Gramvector}).  Then $A\trans{A}$ is diagonal with integer
  entries.
\end{proposition}

Polynomial optimization problems that we consider in the next two
sections are difficult to tackle directly but admit standard SOS
relaxations involving constraints sets \eqref{eq:Gramvector}. In
practice, an SOS relaxation approach boils down to solving linear
semidefinite problems of the form
\begin{equation}\label{eq:pbsdp}
\accol{\min \quad \trans{c}{x}\\
 \quad Ax=b, \ x\in \sdp}
\end{equation}
where and $c\in \RR^{n^2}$, $b\in \RR^m$, $A\in\RR^{m \times n^2}$,
and vector $x$ collects entries of matrix $X \in {\mathbb R}^{n\times
  n}$.  For solving problem \eqref{eq:pbsdp}, we use the three
following the public-domain Matlab implementations:
\begin{enumerate}
\item SeDuMi1.3 implementing the primal-dual interior-point
  algorithm of \cite{sturm-1999}\\ 
  (available on {\tt sedumi.ie.lehigh.edu})
\item MPRW a version of the basic regularization method of
  \cite[Sec.\,5.1]{malick-povh-rendl-wiegele-2009}\\
  (available on {\tt
    www.math.uni-klu.ac.at/or/Software/mprw2.m})
\item SDPNAL0.1 the regularization method of
  \cite{zhao-sun-toh-2010}\\
  (available on {\tt www.math.nus.edu.sg/$\sim$mattohkc/\sdpnal.html})
\end{enumerate}

\noindent
Our goal here is just to show that the regularization methods are
interesting in this context. We simply use default parameter tunings
of the algorithms, with no attempt to adapt them to each particular
problem instances contrary to in \cite{nie-2009}.  With {\tt K.s=n}, the
calling sequences of the three Matlab functions {\tt \sedumi}, {\tt
  \mprw} and {\tt \sdpnal} for solving \eqref{eq:pbsdp} are thus as
follows:
\begin{verbatim}
pars = []; pars.tol = 1e-9;
[x,y] = sedumi(A,b,c,K,pars);
X = reshape(x,K.s,K.s);

tol = 1e-9; C = reshape(c,K.s,K.s);
[X,y] = mprw(A,b,C,1e6,1,tol);

opts = []; opts.tol = 1e-9;
[blk,At,C,B] = read_sedumi(A,b,c,K);
[obj,X,y] = sdpnal(blk,At,C,B,opts);
X = X{1};
\end{verbatim}

Experiments are carried out with Matlab 7.7
running on a Linux PC with Intel Xeon CPU W3520 2.67Ghz using 64 bit
arithmetic and 8GB RAM. Computation times are given in seconds,
with two significant digits only (since our objective is not
a comprehensive accurate benchmarking of codes).

Similarly to \cite{nie-2009}, we will see that, due to lower memory
requirements, regularization methods can solve larger polynomial
optimization problems than classical interior-point methods with the
above setting.

A last note about tolerance. The tolerance parameters {\tt tol} for
the three solvers are set to $10^{-9}$ for all the numerical
experiments (except otherwise stated). Notice though that the meaning
of the tolerance parameter is not the same for two types of algorithms. With
regularization methods, the relative accuracy measured in terms of
primal and residuals (remember \eqref{eq:residual}) is easily
controlled. We stress that lower requirements on the relative accuracy
could result in a significant saving of computational time, and this
could be useful when solving approximately large-scale problems with
\mprw and \sdpnal (see some examples in \cite{nie-2009}).  In
contrast, we observe (as expected) that the iteration count of \sedumi
does not depend significantly on the expected accuracy, measured in
terms of duality gap. Most of the computational time is spent to
generate an approximately feasible primal-dual pair with relatively
small duality gap, and only a few more iterations are required to
refine the accuracy below $10^{-9}$.

\subsection{Testing positivity of polynomials}\label{sec:SOSfeas}

We focus in this section on the very first problem of polynomial
optimization: we would like to know whether a polynomial
\eqref{eq:poly} is positive
\begin{equation}\label{eq:pospoly}
p(v) \geq 0, \quad \text{for all } v \in {\mathbb R}^N.
\end{equation}
In general this is a difficult problem for which no polynomial-time
algorithm is known. It can be relaxed to the easier problem of testing
if $p$ could be expressed as an SOS \eqref{eq:SOSpoly}.
Whenever it holds, then obviously condition (\ref{eq:pospoly})
is satisfied.  The converse is not true in general if $N\geq 2$ and
$d\geq3$, and there are explicit counter-examples; the simplest of
them (the Motzkin polynomial) is studied below.

\subsubsection{Random full-rank polynomial SOS problems}\label{sec:sos}

We consider random polynomial SOS problems
which are constructed so that there is a full-rank orthogonal
Gram matrix $X$ (an interior point) solving problem
(\ref{eq:SOSGram}). We use GloptiPoly 3 (see \cite{lasserre-2009})
to generate matrix $A$ and vector $b$ as follows:
\begin{verbatim}
N = 5; % number of variables
d = 3; % half degree
mpol('v',N,1); % variables
P = msdp(min((v'*v)^d); % construct A matrix
[A,b,c,K] = msedumi(P); % retrieve A and K in SeDuMi format
A = [c';-A]; % constant term and sign change
c = zeros(size(A,2),1); % no objective function
X = orth(randn(K.s)); % random Gram matrix
b = A*X(:); % corresponding right handside vector
\end{verbatim}

On Table \ref{tab:sos} we report execution times (in seconds) required by
\sedumi, \mprw and \sdpnal to solve problem (\ref{eq:Gramvector}) for
$d=3$ (degree six polynomials) and $N=5,\ldots,12$. We also indicate
the size $n$ of matrix $X$ and the number $m$ of constraints (row
dimension of matrix $A$).  We observe that \sedumi is largely
outperformed by \mprw and \sdpnal. We also observe that \mprw is about 4
times slower than \sdpnal, but this is not surprising as \mprw is a
simple prototype (without comments and printing instructions it is
about 50 lines of interpreted Matlab), whereas \sdpnal is a much more
sophisticated package heavily relying on the efficient data handling
and numerical linear algebra routines of the SDPT3 package. We also
recall that \sdpnal makes several iterations of \mprw as preprocessing.

\begin{table}[h]
\centering
\begin{tabular}{ccc|ccc}
$N$ & $n$ & $m$ & SeDuMi & \mprw & SDPNAL \\ \hline
5 & 56 & 462 & 0.29 & 0.03 & 0.05 \\
6 & 84 & 924 & 0.92 & 0.05 & 0.07 \\
7 & 120 & 1716 & 4.8 & 0.13 & 0.10 \\
8 & 165 & 3003 & 25 & 0.35 & 0.16 \\
9 & 220 & 5005 & 110 & 0.66 & 0.25 \\
10 & 286 & 8008 & 410 & 1.3 & 0.43 \\
11 & 364 & 12376 & 1500 & 3.0 & 0.73 \\
12 & 455 & 18564 & $>3600$ & 5.0 & 1.3
\end{tabular}
\caption{Comparative execution times for SOS
problems.\label{tab:sos}}
\end{table}

\subsubsection{Random low-rank polynomial SOS problems}\label{sec:sosrankone}

We consider random polynomial SOS problems which are constructed so
that there is a rank-one Gram matrix $X$ solving problem
(\ref{eq:Gramvector}). For such problems, it is unlikely that there is an
interior point $x$ solving problem (\ref{eq:Gramvector}), and indeed
\sedumi does not find a full-rank solution.  We use the same code as
above, replacing the instruction {\tt X = orth(randn(K.s));} with the
instructions
\begin{verbatim}
X = orth(randn(K.s,1));
X = X*X';
\end{verbatim}

We report execution times (in seconds) in Table \ref{tab:sosrankone}.
In comparison with the problems with interior points of
Table\;\ref{tab:sos}, we observe that all the solvers experience
convergence issues. However, there is still a considerable improvement
in terms of efficiency brought by regularization methods, though
Slater's qualification constraint cannot be invoked to guarantee
convergence.

\begin{table}[h]
\centering
\begin{tabular}{ccc|ccc}
$N$ & $n$ & $m$ & SeDuMi & MPRW & SDPNAL \\ \hline
5 & 56 & 462 & 0.83 & 0.20 & 0.21 \\
6 & 84 & 924 & 0.85 & 0.32 & 0.28 \\
7 & 120 & 1716 & 16 & 2.1 & 0.51 \\
8 & 165 & 3003 & 61 & 4.8 & 0.98 \\
9 & 220 & 5005 & 330 & 12 & 1.2 \\
10 & 286 & 8008 & 1300 & 24 & 2.5 \\
11 & 364 & 12376 & $>3600$ & 50 & 3.5 \\
12 & 455 & 18564 & $>3600$ & 110 & 6.6
\end{tabular}
\caption{Comparative execution times for low-rank SOS
problems.\label{tab:sosrankone}}
\end{table}

\subsubsection{Motzkin's polynomial}\label{sec:motzkin}

We study a well-known bivariate ($N=2$) polynomial of sixth degree
($d=3$) which is non-negative but
cannot be written as a polynomial SOS, namely
Motzkin's polynomial
\[
p_0(v) = 1+v^2_1v^2_2(v^2_1+v^2_2-3)
\]
see \cite{laurent-2009} or \cite{lasserre-2009}.
This polynomial achieves its minimum zero
at the four points $v_1=\pm 1$, $v_2=\pm 1$.
In a basis of monomials of degree up to $3$
there is no Gram matrix $X$
solving problem (\ref{eq:Gramvector}).
However, it was observed in \cite{henrion-lasserre-2005}
and later on shown theoretically in \cite{lasserre-2006}
that the perturbed polynomial
\[
p_0(v) + \varepsilon p_1(v)
\]
can be represented as a polynomial SOS (with full-rank Gram matrix)
provided the degree of the perturbation polynomial $p_1(v)$ is high
enough, inversely proportional to scalar $\varepsilon>0$.  In some
sense, this can be interpreted as a regularization procedure as in
\cite{henrion-malick-2009}.  Practically speaking, since semidefinite
programming solvers use inexact operations (floating point
arithmetic), it is not necessary to perturb explicitly the data. It is
enough to choose a basis $\pi(v)$ of sufficiently high degree $d>3$ in
relation (\ref{eq:SOSGram}), and higher-order perturbations are
automatically introduced by the algorithm.

We use the following GloptiPoly3 instructions
to generate data $A$, $b$ for increasing values of $d$:
\begin{verbatim}
d = 8; % half degree
mpol v1 v2
p = 1+v1^2*v2^2*(v1^2+v2^2-3);
P = msdp(min(p),d);
[A,b,c,K,b0] = msedumi(P);
A = [c';-A];
b = [-b0;-b];
c = zeros(size(A,2),1);
\end{verbatim}

For this problem, we set {\tt tol=1e-6} for the three solvers.
When $d=3,4,5,6$, \sedumi takes less than 0.1 seconds to
detect that problem (\ref{eq:Gramvector}) is infeasible,
and it provides a Farkas dual certificate vector $y \in -{\mathcal K}$
such that $\trans{b} y = 1$.
When $d=7,8,9,10$, \sedumi takes less than 0.5 seconds
to return a vector $x$ such that the primal residual
$\|Ax-b\|_2/\|b\|_2$ is less than $10^{-9}$
and the dual objective function $\trans{b} y$ is less than
$10^{-9}$ in absolute value.

The behavior of \sdpnal and \mprw is more erratic, and convergence
issues occur, as shown by the execution times (in seconds) of
Table \ref{tab:motzkin}.
For $d=3,4,5$, \mprw stops after $10^6$ iterations, as there is no mechanism
to detect infeasibility in this prototype software.

\begin{table}[h]
\centering
\begin{tabular}{c|ccc}
$d$ & time & $\|Ax-b\|_2/\|b\|_2$ & $\trans{b}y$ \\ \hline
3 & - & - & - \\
4 & - & - & - \\
5 & - & - & - \\
6 & 15 & $6.20 \cdot 10^{-6}$ & $1.12 \cdot 10^{-6}$ \\
7 & 25 & $6.03 \cdot 10^{-7}$ & $6.81 \cdot 10^{-7}$ \\
8 & 26 & $5.80 \cdot 10^{-6}$ & $-4.08 \cdot 10^{-7}$ \\
9 & 34 & $1.01 \cdot 10^{-6}$ & $-1.45 \cdot 10^{-7}$ \\
10 & 75 & $5.42 \cdot 10^{-7}$ & $-1.58 \cdot 10^{-7}$
\end{tabular}
\qquad\quad
\begin{tabular}{c|ccc}
$d$ & time & $\|Ax-b\|_2/\|b\|_2$ & $\trans{b}y$ \\ \hline
3 & 5.1 & $4.28 \cdot 10^{-3}$ & 33.4 \\
4 & 9.2 & $1.56 \cdot 10^{-4}$ & 0.832 \\
5 & 3.5 & $4.59 \cdot 10^{-6}$ & $4.37 \cdot 10^{-5}$ \\
6 & 4.6 & $6.33 \cdot 10^{-6}$ & $1.05 \cdot 10^{-6}$ \\
7 & 5.7 & $8.95 \cdot 10^{-6}$ & $3.86 \cdot 10^{-7}$ \\
8 & 5.9 & $2.79 \cdot 10^{-6}$ & $-3.46 \cdot 10^{-7}$ \\
9 & 7.9 & $2.54 \cdot 10^{-6}$ & $-3.25 \cdot 10^{-7}$ \\
10 & 8.8 & $1.88 \cdot 10^{-6}$ & $-1.34 \cdot 10^{-7}$
\end{tabular}
\caption{Behavior of \mprw (left) and \sdpnal (right) for Motzkin's polynomial.\label{tab:motzkin}}
\end{table}

\subsubsection{Regularization vs projection}

Though it solves linear semidefinite problems, using regularization
techniques somehow generalizes and enhances the idea
\cite{henrion-malick-2009} to using projection methods directly for
SOS feasibility problems. With this approach indeed, a question is to
find a good point to project; taking systematically the zero matrix
gives interesting results but could be greatly
enhanced. Regularization methods provide a numerical solution to this:
doing a sequence of truncated projections allows to keep the
projection idea while getting rid of the question of the initial point
to project. The behaviour of SDPNAL is interesting with this respect:
it does first a preprocessing of several alternating direction
iterations to get a meaningful point, then follows by projection-like
iterations. In practice, we observe usually a very few iterations, and
often one.  For example, to decide whether the (admittedly trivial)
polynomial $p(v)=\sum_{i=1}^{10}v_i^{10}$ is SOS, the SDP problem
dimensions are $n=3003$ and $m=184756$, and after 90 seconds and only
one projection-like iteration, \sdpnal provides a vector $x$
satisfying $\|Ax-b\|_2/\|b\|_2 \approx 1.4\cdot 10^{-10}$.

\subsection{Unconstrained polynomial minimization}\label{sec:SOSmin}

In this section we study global minimization problems
\begin{equation}\label{eq:minpoly}
p^* = \min_{v \in {\mathbb R}^N} p(v)
\end{equation}
where $p(v)$ is a given polynomial. For this problem, a semidefinite
relaxation readily follows from the observation that
\[
\begin{array}{rcll}
p^* & = & \max_{\underline{p}} & \underline{p} \\
& & \mathrm{s.t.} & p(v)-\underline{p} \geq 0, \quad \forall v \in {\mathbb R}^N
\end{array}
\]
and by relaxing the above non-negativity constraint by the semidefinite
programming constraint that polynomial $p(v)-\underline{p}$ is SOS,
see \cite{laurent-2009} and \cite{lasserre-2009}.

\subsubsection{Random polynomial minimization problems}

We generate well-behaved instances of unconstrained polynomial minimization
problems (\ref{eq:minpoly}) with
\[
p(v) = p_0(v) + \sum_{i=1}^N v_i^{2d}
\]
where $p_0(v)$ is a random polynomial of total degree strictly less than $2d$.
The leading term $\sum_{i=1}^N v_i^{2d}$
ensures coercivity of $p(v)$ and hence existence of a global minimum
in (\ref{eq:minpoly}).
We use the following GloptiPoly 3 script to generate our examples:
\begin{verbatim}
N = 10;
mpol('v',N,1);
b = mmon(v,0,2*d-1); % degree up to 2d-1
p0 = randn(1,length(b)); p0 = p0/norm(p0);
p = p0*b + sum(mmon(v,d).^2);
P = msdp(min(p));
[A,b,c,K] = msedumi(P);
\end{verbatim}

In Table \ref{tab:minpol} we report comparative execution times (in
seconds) for $d=2$ and various values of $N$, for solving the
semidefinite relaxation. It turns out that for these
generic problems, we observe that global optimality is always
certified with a rank-one moment matrix \cite{henrion-lasserre-2005}.
Both \mprw and \sdpnal largely outperform \sedumi on these examples.

\begin{table}[h]
\centering
\begin{tabular}{ccc|ccc}
$N$ & $n$ & $m$ & SeDuMi & MPRW & SDPNAL \\ \hline
5 & 21 & 126 & 0.09 & 0.05 & 0.18 \\
6 & 28 & 209 & 0.11 & 0.07 & 0.18 \\
7 & 36 & 329 & 0.24 & 0.12 & 0.20 \\
8 & 45 & 494 & 0.36 & 0.19 & 0.22 \\
9 & 55 & 714 & 0.77 & 0.28 & 0.26 \\
10 & 66 & 1000 & 1.9 & 0.45 & 0.29 \\
11 & 78 & 1364 & 5.0 & 0.78 & 0.36 \\
12 & 91 & 1819 & 11 & 1.1 & 0.41 \\
13 & 105 & 2379 & 20 & 1.6 & 0.47 \\
14 & 120 & 3059 & 42 & 2.3 & 0.65 \\
15 & 136 & 3875 & 74 & 3.0 & 0.68
\end{tabular}
\caption{Comparative execution times for semidefinite relaxations of
random polynomial minimization problems.\label{tab:minpol}}
\end{table}

\subsubsection{A structured example}\label{sec:nie35}

Consider the problem studied in \cite[Example 3.5]{nie-2009}, that is
(\ref{eq:minpoly}) with
\[
p(v) = \sum_{i=1}^N \left(1-\sum_{j=1}^i
(v_j+v_j^2)\right)^2+\left(1-\sum_{j=1}^N(v_j+v_j^3)\right)^2.
\]
We solve the semidefinite relaxation for
increasing values of $N$.  We collect comparative execution times on
Table \ref{tab:nie35} for this example. 
For example, when $N=10$, \sdpnal resp. \mprw returns a point $x$
such that $\|Ax-b\|/\|b\|$ is equal to $1.4\cdot 10^{-9}$ resp.
$2.6\cdot 10^{-10}$
and the minimum eigenvalue of $X$ is equal to zero to machine precision.

We observe again a considerable improvement in terms of performance brought by
regularization methods in comparison with a classical interior-point
method.  For larger instances, most of the computation time of \sedumi
is spent for memory swapping when constructing and handling large
matrices.  We refer to the recent extensive numerical work of
\cite{nie-2009} for various structured problems.

\begin{table}[h]
\centering
\begin{tabular}{ccc|ccc}
$N$ & $n$ & $m$ & SeDuMi & MPRW & SDPNAL \\ \hline
5 & 56 & 461 & 0.50 & 0.54 & 0.60 \\
6 & 84 & 923 & 2.5 & 1.2 & 1.2 \\
7 & 120 & 1715 & 14 & 5.0 & 2.6 \\
8 & 165 & 3002 & 92 & 19 & 6.7 \\
9 & 220 & 5004 & 410 & 65 & 22 \\
10 & 286 & 8007 & 1800 & 200 & 71 \\
11 & 364 & 12375 & 7162 & 490 & 150 \\
12 & 455 & 18563 & $>7200$ & 1500 & 530 \\
13 & 560 & 27131 & $>7200$ & 3500 & 2300 \\
14 & 680 & 38760 & $>7200$ & $>7200$ & 9900 \\
\end{tabular}
\caption{Comparative execution times for semidefinite relaxations of
a larger polynomial minimization problem.\label{tab:nie35}}
\end{table}

\section*{Acknowledgements}

The work of the first author was partly supported by
Research Programme MSM6840770038 of the Czech Ministry of Education
and by Project 103/10/0628 of the Grant Agency of the Czech Republic.

\end{document}